\documentclass[10pt,a4paper]{article} 
\usepackage[latin1]{inputenc}
\usepackage[T1]{fontenc}
\usepackage{lmodern}
\usepackage{indentfirst}
\usepackage[francais,english]{babel}
\usepackage{amsmath}
\usepackage{amsthm}
\usepackage{amsfonts}
\usepackage[vmargin=2.5cm,hmargin=3.5cm]{geometry}
\usepackage{amssymb}
\usepackage{pstricks-add}
\usepackage{graphicx} 
\usepackage{color}
\usepackage{thmbox}
\usepackage{wallpaper}
\usepackage[french]{minitoc}
\usepackage{fancyhdr} 
\color{black}
\usepackage{fancyhdr}
\usepackage{url}
\usepackage{hyperref} 
\hypersetup{
colorlinks=true, 
breaklinks=true, 
urlcolor= blue, 
linkcolor= red, 
citecolor=black, 
pdftitle={}, 
pdfauthor={Lahoucine Elaissaoui}, 
pdfsubject={Simulation} 
} 
\def\R {\mathbb{R}}
\def\N {\mathbb{N}}

\def\C {\mathbb{C}}
\def\L {\mathcal{L}}

\def\vb {\mathcal{V}_b\mathcal{C}_g}

\def\eps {\varepsilon}

\def\P {\mathcal{P}}

\renewcommand{\abstractname}{\sc{Abstract}}

\newcommand{\equa}[1]{\begin{equation} #1 \end{equation}}

\newcommand{\theo}[2]{\begin{thee}\fbox{\begin{minipage}{\textwidth} \par #1 \end{minipage}} \end{thee} \vspace{0.1cm} \textbf{\textit{\sc{Preuve:}}} \vspace{0.5cm} \par  #2  \begin{flushright}$\blacksquare$\end{flushright} }

\newcommand{\lem}[2]{\begin{lemm}\fbox{\begin{minipage}{\textwidth} \par #1 \end{minipage}} \end{lemm} \vspace{0.1cm} \textbf{\textit{\sc{preuve:}}} \vspace{0.5cm} \par  #2  \begin{flushright}$\blacksquare$\end{flushright} }

\newcommand{\corol}[2]{\begin{coro}\fbox{\begin{minipage}{\textwidth} \par #1 \end{minipage}} \end{coro} \vspace{0.1cm} \textbf{\textit{\sc{preuve:}}} \vspace{0.5cm} \par  #2  \begin{flushright}$\blacksquare$\end{flushright} }

\newtheorem{thee}{\sc Théorème}[section]
\newtheorem{propo}{Proposition}[section]
\newtheorem{lemm}{\sc Lemme}[section]
\newtheorem{defin}{Définition}[section]

\newtheorem{remr}{Remarque}[section]
\newtheorem{propr}{Propriétés}[section]
\newtheorem{coro}{\sc Corollaire}[section]
\title{  \textsc{Tauberian Theorem of Laplace Transformation \\ And \\ Application of Prime Number Theorem }  }
\author{\color{blue}{\textsc{Lahoucine Elaissaoui}} \\ \texttt{lahoucine.elaissaoui@stud-mail.uni-wuerzburg.de} \\ \texttt{lahoumaths@gmail.com} }
\begin{document}
\maketitle

\vspace*{3cm}
\renewcommand{\abstractname}{\sc{Résumé}}
 \begin{abstract}

 {\footnotesize Dans cet article je donnerai une nouvelle démonstration courte et directe pour le Théorème des Nombres Premiers. C'est vrai que ce théorème a été complétement démontré au début du 20ème siecle mais la démonstration était basé sur des résultats élémentaires (théorème de \textbf{Chebyshev}) et aussi analytiques compliqués (théorème de \textbf{Ikehrara}), mais ici j'ai pas utilisé le théorème de Chebyshev ainsi que j'ai remplacé et j'ai généralisé le théorème de \textbf{Ikehara} grâce à la notion des fonctions à variation bornée qui est ancienne mais récent dans la théorie analytique des nombres.}
  
\end{abstract}

\begin{center}
\section{\sc Préliminaires}
\end{center}
\subsection{Les fonctions à variation bornée}
Les fonctions à variation bornée joue un rôle très important dans la théorie de l'intégration au sens de \textbf{Stieltjes}, ici on va s'interésser à les fonctions à variation bornée sur $\R^+$ à valeurs complexes. Soit $x$ un réel positif et soit $(x_k)_{k=0,\cdots, n}$ une suite finie et strictement croissante des réels de l'intervalle $[0,x]$ tels que $0=x_0 < x_1<x_2< \cdots < x_n = x$ est une subdivision de l'intervalle $[0,x]$, on note $\Sigma$ pour cette subdivion et $\mathcal{S}([0,x])$ l'ensemble de toutes les subdivions possibles de $[0,x]$.  \textit{La fonction variation totale} d'une fonction complexe définie sur $\R^+$, notée $T_f$, est la fonction définie par \equa{T_f(x) := \sup_{\Sigma \in \mathcal{S}([0,x])} \sum_{k=1}^n|f(x_k) - f(x_{k-1})|\label{e1}} 

Il est bien clair que la fonction $T_f$ est une fonction croissante sur $\R^+$, par conséquent si $T_f$ est majorée sur $\R^+$ alors on dira que $f$ est à \textit{variation bornée} sur $\R^+$ et on note $$V(f) = \lim_{x\to + \infty} T_f(x)  \in \R^+$$ pour \textit{la variation totale} de la fonction $f$.

\begin{propr}
\begin{itemize}
\item Toute fonction $g$ de classe $\mathcal{C}^1$ sur $\R^+$ à valeurs complexes telle que $g' \in L^1(\R^+)$ est à variation bornée, en effet, pour une subdivision $\Sigma: 0=x_0 < x_1<\cdots < x_n=x$ et puisque $g$ est continue sur chaque intervalle $[x_{i-1},x_i]$ (pour $i=1,\cdots n$) et dérivable sur leurs interieurs topologique alors d'après le théorème des accroissements finis il existe des $c_i$ dans $]x_{i-1},x_i[$ tels que $$|g(x_i) - g(x_{i-1})|=|g'(c_i)||x_{i} - x_{i-1}| $$

D'où $$T_g(x) = \sup_{\Sigma \in \mathcal{S}([0,x])} \sum_{k=1}^n|g'(c_i)||x_{i} - x_{i-1}|$$

Or cette somme est une somme de \textbf{Darboux} ce qu'on peut déduire grâce à l'intégrale de \textbf{Riemann} que $$T_g(x) = \int_0^x|g'(t)| dt $$

Donc $$V(g) = \int_{0}^{+\infty}|g'(t)| dt $$
qui est finie puisque $g' \in L^1(\R^+)$, alors $g$ est à variation bornée sur $\R^+$.
 \item toute fonction à variation bornée sur $\R^+$ est bornée sur $\R^+$, en effet, soit $f$ une fonction à variation bornée sur $\R^+$ alors pour un $x\geq 0$ \begin{align*}|f(x) - f(0)| &= \left|\sum_{k=1}^n f(x_i) - f(x_{i-1})\right| \\ &\leq \sum_{k=1}^n|f(x_i) - f(x_{i-1})| \\ &\leq T_f(x) \\ &\leq V(f) <+\infty  \end{align*}
 
 Alors $f$ est bornée sur $\R^+$.
\end{itemize} \label{prop1}\end{propr}

On dit qu'une fonction $f$ définie de $\R^+$ à valeurs complexes admet une limite à gauche en $x \in \R^+$, notée $f(x^-)$ si à tout $\eps >0$ on peut associer un $0 \leq \delta < x$ tel que  
$$a < t < x \Longrightarrow |f(t) - f(x^-)| < \eps $$

Et en plus si $f(x^-) = f(x)$ on dit que $f$ est continue à gauche en $x$. 

On note pour $\vb$ la classe des fonctions, définies de $\R^+$ à valeurs complexes, à variation bornée, continues à gauche en tout point de $\R^+$ et qui s'annullent en $0$.

\subsection{Intégrale de Lebesgue-Stieltjes}

Le théorème 8.14 page 156 du livre [Rud] a établi le lien entre la théorie de la mesure et la théorie des fonctions à variation bornée. Donc d'après le même théorème, soit $f \in \vb$ alors il existe une unique mesure complexe de Borel $\mu_f$ telle que \equa{f(x) = \mu_f([0,x[), \qquad \forall x\geq 0 \label{e2}}

Et en plus pour tout $x \in \R^+$ on a \equa{T_f(x)=|\mu_f|([0,x[)}

Où $|\mu_f|$ est une mesure positive de Borel, dite \textit{la variation totale de la mesure complexe $\mu_f$}, qui est \textit{finie} d'après le théorème 6.4 page 114 de [2].
\begin{remr}
\begin{itemize}
\item On peut facilement montrer que $|\mu_f|$ est finie autant que $f\in \vb$, en effet, soit $f\in \vb$ alors \begin{align*} |\mu_f|(\R^+) &= \lim_{x\to +\infty} |\mu_f|([0,x[) \\ &= \lim_{x\to +\infty}T_f(x) \\&= V(f) <+ \infty\end{align*}

\item D'une autre part, si $f$ est à valeurs dans $\R$ alors $\mu_f$ est dite \textit{une mesure signée} alors de la même manière on démontre que cette mesure est finie.
 \item Soit $f\in \vb$, si $y>x$ alors \begin{align*}f(y) - f(x) &= \mu_f([0,y[) - \mu_f([0,x[)\\ &= \mu_f([x,y[) \end{align*}
 
 Donc $$\mu_f(\{x\}) = f(x^+) - f(x) $$
 
 D'où  $f$ est continue en $x$ si et seulement si $$\mu_f(\{x\}) = 0 $$
\end{itemize}
\end{remr}
Le théorème de \textbf{Radon-Nikodym}, voir le théorème 6.12 page 120 de [2], assure que pour toute mesure complexe $\mu$ il existe une fonction mesurable complexe $h$ de module égal à $1$ telle que $$d\mu = h d|\mu| .$$

 Ainsi, on déduit que pour toute fonction $g:\R^+ \longrightarrow \C$ mesurable et bornée sur $\R^+$ on a $g \in L_{\mu_f}^1(\R^+)$ où $f\in \vb$. En effet: \begin{align*} \left| \int_{\R^+} g d\mu_f \right| &\leq \int_{\R^+} |g| d|\mu_f| \\ &\leq \|g\|_{\infty} |\mu_f|(\R^+)\\ &< +\infty \end{align*}
 
 Où $$\|g\|_{\infty} = \sup_{x \in \R^+}|g(x)|.$$
 
 Maintenant, d'après le théorème 6.1.4 du livre [1] on constate que pour $f \in \vb$ on a\equa{\int_0^x df(t) = \mu_f([0,x[), \qquad x\geq 0 \label{e4}}
 
 Soient donc $f\in \vb$ et $g:\R^+ \longrightarrow \C$ une fonction de classe $\mathcal{C}^1$ sur $\R^+$ telle que $g' \in L^1(\R^+)$, alors d'après le théorème 6.2.2 (grâce au résultat \ref{e4}) du même livre on démontre que \equa{\int_0^{+\infty} g(t) df(t) = \mu_{fg}(\R^+) - \int_0^{+\infty} f(t)g'(t) dt \label{e5}}
 
 la mesure complexe $\mu_{fg}$ a bien un sens, en effet d'après les propriétés \ref{prop1} on démontre que $g$ est à variation bornée or le produit de deux éléments de $\vb$ est un élément de $\vb$ alors  $fg \in \vb$ (car $fg$ est à variation bornée et continue à gauche à chaque point de $\R^+$ et $(fg)(0) =0$), et en plus $$\mu_{fg}([0,x[) = f(x)g(x), \qquad \forall x\in \R^+ $$
 
 Et $$|\mu_{fg}(\R^+)| \leq |\mu_{fg}|(\R^+) = \lim_{x\to + \infty} T_{fg}(x)< + \infty $$
 
 \section{Théorème Tauberien de la transformation de Laplace complexe}
 Dans tout ce qui suit $s = \sigma + it$ où $\sigma , t \in \R$ est un nombre complexe et $\rho$ est une fonction de la classe $\vb^* := \{f \in \vb , \quad  \ \Im(f) = 0 \}$. Ainsi $\C_*^+$ est l'ensemble des nombres complexes de partie réelle strictement positive.
 
 On définit la transformation de Laplace-Stieltjes de la fonction $\rho$ par $$\L_{\rho}^*(s) = \int_0^{+\infty} e^{-sx}d\rho(x), \qquad \sigma > 0 $$
 
 il est bien clair d'après ce qui précéde, puisque $x\mapsto e^{-sx}$ est continue et bornée sur $\R^+$ pour tout $\sigma > 0$, que la fonction $\L_{\rho}^*$ est bien définie.
 
 \lem{Soit $\rho \in \vb^*$ alors $$\lim_{s \to 0^+} \L_{\rho}^*(s) = \lim_{+ \infty} \rho .$$\label{L1}}{On pose pour tout $(x,s) \in \R^+ \times \C_+^*$ $$\phi(x,s) = e^{-sx} $$
 
 Alors on a
 
 \begin{itemize}
\item Pour tout $x\geq 0$ la fonction $s\mapsto \phi(x,s)$ est continue en $0^+$.
\item Pour tout $\sigma > 0$ la fonction $x \mapsto \phi(x,s)$ est continue donc mésurable sur $\R^+$.
\item Pour tout $\sigma > 0$ et pour $d\rho$-presque tout $x \in \R^+$ on a $$|\phi(x,s)|\leq 1$$

Où $1 \in L_{d\rho}^1(\R^+)$ car $$\int_{\R^+} d\rho(x) = \mu_{\rho}(\R^+) <+ \infty $$

Alors la fonction $x \mapsto \phi(x,s)$ est $d\rho$-intégrable sur $\R^+$ et la fonction $\L_{\rho}^*$ est est continue en $0^+$. Donc \begin{align*} \lim_{s \to 0^+} \L_{\rho}^*(s) &= \L_{\rho}^*(0) \\ &= \int_{\R^+}d\rho (x) \\ &= \lim_{x\to + \infty} \mu_{\rho}([0,x[) \qquad \text{d'après \ref{e4}}  \\ &= \lim_{+ \infty} \rho \qquad \qquad \qquad \quad \! \! \text{d'après \ref{e2}} \end{align*}
\end{itemize}}

Maintenant on définit la transformation de Laplace complexe d'une fonction $\rho \in \vb^*$ par $$\L_{\rho}(s) = \int_{0}^{+ \infty} \rho(x)e^{-sx}dx, \qquad \sigma > 0 .$$

La fonction $\L_{\rho}$ est bien définie, en effet puisque $\rho \in \vb^*$ alors $\rho$ est bornée sur $\R^+$ en plus $$\left|\int_0^{+ \infty} \rho(x)e^{-sx}dx \right| \leq \| \rho\|_{\infty} \int_0^{+ \infty} e^{- \sigma x} dx = \frac{\|\rho\|_{\infty}}{\sigma} < +\infty$$ 

\theo{Soit $\rho \in \vb^*$, on suppose que $\L_{\rho}$ est holomorphe sur $\{\sigma > 0 \}$ et admet un prolongement analytique sur $\{\sigma \geq 0 \}$ avec un pôle simple au point $0$. Alors on a $$\lim_{+ \infty}\rho = Res(\L_{\rho},0) .$$\label{T1}}{Soit $s\in \C_+^*$, d'après l'équation \ref{e5} on a $$\L_{\rho}^*(s) = \mu_{\rho e^{-s \cdot}}(\R^+) + s \int_0^{+\infty} \rho(x) e^{-sx} dx$$

Or $$|\mu_{\rho e^{-s \cdot}}(\R^+)| = |\mu_{\rho}(\R^+)|\lim_{x \to + \infty} e^{-\sigma x} = 0  $$

Donc $$\L_{\rho}^*(s) = s \L_{\rho}(s) $$

Passons à la limite $s \to 0^+$ on a d'après le \textsc{Lemme} \ref{L1} $$\lim_{x\to + \infty} \rho(x) =Res(\L_{\rho},0)  $$ }

D'une manière générale, soit $\alpha$ un réel positif alors il est clair, d'après ce qui précéde, que pour tout $\rho \in \vb^*$ on a $\varrho(x) = \rho(x) e^{-\alpha x}$ est un élément de $\vb^*$. Ainsi on déduit le résultat suivant:

\corol{Soient $\alpha \in \R^+$ et $\rho \in \vb^*$. On suppose que la fonction $\L_{\rho}$ est holomorphe sur $\{\sigma > \alpha \}$ et admet un prolongement analytique sur $\{\sigma \geq \alpha \}$ avec un seul pôle simple en $s=\alpha$ alors on a $$\rho (x) \underset{x \to + \infty}{\sim} Res(\L_{\rho}, \alpha) e^{\alpha x}. $$ \label{C1}}{Soit $\sigma > \alpha$, alors $$\L_{\rho}(s) = \int_0^{+\infty} \rho(x)e^{-sx}dx = \int_{0}^{+\infty} \rho(x) e^{-\alpha x} e^{-(s-\alpha)x}dx $$

On pose $$\varrho(x) = \rho(x)e^{- \alpha x}, \qquad \forall x \geq 0  .$$

Alors $$\L_{\rho}(s) = \int_0^{+ \infty} \varrho(x)e^{-(s-\alpha)x}dx = \L_{\varrho}(s - \alpha) $$

Donc $$(s-\alpha) \L_{\rho}(s) = (s-\alpha)\L_{\varrho}(s- \alpha) = z \L_{\varrho}(z) $$

D'où quand $s \to \alpha$ on aura $z \to 0$ et d'après le \textsc{Théorème} \ref{T1} on a $$\lim_{x \to + \infty} \varrho(x) = Res(\L_{\varrho}(z),z=0) = Res(\L_{\rho},\alpha) $$

Alors $$\rho(x) \underset{x \to + \infty}{\sim} Res(\L_{\rho},\alpha) e^{\alpha x} .$$ }

\section{Théorème des Nombres Premiers (nouvelle démonstration)}

Soit $\chi: \N^* \longrightarrow \R^+$ une fonction arithmétique positive, on pose pour tout $x \in (1,+\infty)$ $$f(x) = \sum_{1 \leq n < x}\chi(n) \qquad \text{et} \qquad f(1) =0 .$$

Il est clair que la fonction $f$ est croissante sur $[1,+\infty)$ et continue à gauche en tout point de $[1,+\infty)$. Ainsi, les points de discontinuité de $f$ sont des éléments de $\N^*$. Si $f$ est continue en $x\in \N^*$ alors on aura $$f(x^+) = f(x) $$

Donc \begin{align*} 0 &= f(x^+) - f(x) \\ &= \sum_{x \leq n < x^+}\chi(n) \\ &= \chi(x) \end{align*}

Alors \equa{f \ \text{est \ continue \ en } \ x \in \N^* \Longleftrightarrow \chi(x) = 0 \label{R1}}

Soit maintenant $(a_k)_{k \in \N}$ une suite croissante des points de discontinuité de la fonction $f$ sur $[1,+\infty)$ alors $f$ est constante sur chaque intervalle $I_k=(a_{k-1},a_{k}]$ (où $k \in \N^*$). En effet, soit $k\in \N^*$ s'il existe $n \in I_k$ tel que $f$ est continue en $n$ alors d'après \ref{R1} $\chi(n) = 0$ ainsi et d'une manière générale soit $(\beta_i)_{i \in \N^*}$ une suite strictement croissante des entiers de $\overset{\circ}{I_k}$ (l'interieur de $I_k$), alors $f$ est continue en chaque $\beta_i$ d'où $\chi(\beta_i)=0$ pour tout $i=1,2,\cdots$ et en conséquent pour tout $x \in (a_{k-1},a_k]$ on a $f(x) = f(a_{k-1}^+)$ ($k \in \N^*$).

Soient $\alpha >1$ un réel et $\rho$ la fonction définie sur $\R^+$ par $$\rho(x) = f\left( e^{x} \right) e^{-\alpha x} .$$

Soit $k$ un entier strictement positif on note $ (\lambda_k)_{k \in \N}$ pour la suite croissante des points de discontinuité de la fonction $\rho$ sur $\R^+$ ($\lambda_k = \log a_k \in \log\N^*$). Alors la fonction $\rho$ est décroissante sur chaque intervalle $J_k = (\lambda_{k-1},\lambda_k]$, en effet: soient $x,y \in J_k$ tels que $x>y$, donc puisque $f$ est constante ($\equiv c_k$) sur $J_k$ alors $\rho(x) - \rho(y) = c_k \left(e^{-\alpha x} - e^{-\alpha y}\right) < 0$ d'où $\rho$ est strictement décroissante sur $J_k$ pour tout $k \in \N^*$. D'une autre part, pour tout $k \in \N^*$ $$\rho(\lambda_k^+) > \rho(\lambda_k). $$

En effet, puisque la fonction $x \mapsto e^{- \alpha x}$ est continue sur $\R^+$ alors $e^{-\alpha \lambda_k^+} = e^{- \alpha \lambda_k}$ donc $$\rho(\lambda_k^+)-\rho(\lambda_k) = (f(a_k^+) - f(a_k))e^{- \alpha \lambda_k} $$
et puisque $f$ est discontinue en $a_k$ et croissante sur $[1,+\infty)$ alors $f(a_k^+) > f(a_k)$. D'où $$\rho(\lambda_k^+) > \rho(\lambda_k). $$

\lem{Soit $\alpha > 1$ alors $$\sum_{n \geq 1}\frac{\chi(n)}{n^{\alpha}} < + \infty \Longrightarrow  \rho \in \vb^*  $$\label{L1}}{Soit $x \in \R^+$, on pose $0=x_0 < x_1 < \cdots < x_n = x$ une subdivision de l'intervalle $[0,x]$ et on note pour $m$ le plus grand entier naturel non nul tel que $\lambda_{m-1} < x \leq \lambda_m$ où les $(\lambda_k)_{k \in \N}$ sont les points de discontinuité de la fonction $\rho$ définis précédamment, alors  $$\sum_{i=1}^n |\rho(x_i) - \rho(x_{i-1})| = \sum_{k=0}^m\sum_{\underset{x_i \in J_k}{i=1}}^n|\rho(x_i) - \rho(x_{i-1})| $$
où $(J_k)_{k \in \N^*}$ sont les intervalles $(\lambda_{k-1},\lambda_k]$ et $J_0=[0,\lambda_0]$, et on note bien que $\displaystyle \cup_{k=0}^mJ_k =[0, \lambda_m] $ donc puisque $\rho$ est strictement décroissante sur chaque $J_k$ alors

 \begin{align*}\sum_{k=0}^m\sum_{\underset{x_i \in J_k}{i=1}}^n|\rho(x_i) - \rho(x_{i-1})| & \leq \rho(0) - \rho(\lambda_0) + \sum_{k=1}^m\left(\rho(\lambda_{k-1}^+) - \rho(\lambda_k)\right) - (\rho(x) - \rho(\lambda_m)) \\ &=-\rho(\lambda_0)+\rho(\lambda_0^+)-\rho(\lambda_1)+\rho(\lambda_1^+)+ \cdots -\rho(\lambda_m) - \rho(x) + \rho(\lambda_m)\\ &= - \rho(x) + \sum_{k=0}^{m-1}\left(\rho(\lambda_k^+) -\rho(\lambda_k)\right) \\ &= - \rho(x) + \sum_{k=0}^{m-1} \frac{f(a_k^+) - f(a_k)}{a_k^{\alpha}} \\ &= -\rho(x) + \sum_{k=0}^{m-1}\frac{\chi(a_k)}{a_k^{\alpha}}\end{align*}

Où les $(a_k)_{k \in \N}$ sont les points de discontinuité de la fonction $f$ et qui sont des éléments de $\N^*$. Donc $$\sum_{k=0}^{m-1}\frac{\chi(a_k)}{a_k^{\alpha}} \leq \sum_{1\leq \ell < e^x} \frac{\chi(\ell)}{\ell^{\alpha}} $$

D'où $$ \sum_{i=1}^n |\rho(x_i) - \rho(x_{i-1})| \leq - \rho(x) + \sum_{1\leq \ell < e^x} \frac{\chi(\ell)}{\ell^{\alpha}}.$$

Alors $$T_{\rho}(x) \leq - \rho(x) + \sum_{1 \leq \ell <e^x}\frac{\chi(\ell)}{\ell^{\alpha}}  .$$

Or puisque $\rho$ est une fonction positive alors $$T_{\rho}(x) \leq \sum_{1 \leq \ell < e^x} \frac{\chi(\ell)}{\ell^{\alpha}} .$$

Donc $$\text{la série} \ \sum_{n \geq 1} \frac{\chi(n)}{n^{\alpha}} \ \text{converge} \ \Longrightarrow \rho \in \vb^* .$$}

Sans perte de généralité le résultat est vrai pour toute fonction arithmétique $\chi:\N^* \longrightarrow \R$ croissante. Dans ce cas, le \textsc{Lemme} \ref{L1} peut être reformulé: $$\text{la série} \ \sum_{n \geq 1} \frac{\chi(n)}{n^{\alpha}} \ \text{est absolument convergente} \Longrightarrow \rho \in \vb^*. $$
où $\alpha > 1$.

Maintenant on arrive au résultat le plus important dans cette section:
\newpage

\theo{Soit $\chi : \N^* \longrightarrow \R^+$ une fonction arithmétique et soit $\alpha> 1$ le plus petit réel tel que la série $\displaystyle \sum_{n \geq 1} \frac{\chi(n)}{n^{\alpha}}$ soit convergente. On pose $\displaystyle f(x) = \sum_{1 \leq n < x} \chi(n)$ pour tout $x \geq 1$ où $f(1)=0$ et on suppose que $\L_{\rho}$ est holomorphe sur le demi-plan complexe $\{\sigma \geq \beta\}$  sauf au seul pôle simple en $s=\beta \geq 0$ alors on a $$f(x) \underset{x \to + \infty}{\sim} Res(\L_{\rho},\beta) x^{\alpha + \beta} .$$
Où $\rho(x) = f(e^x)e^{-\alpha x}$ pour tout $x \in \R^+$. \label{T2}}{Soit $\alpha > 1$ un réel tel que la série du terme générale $\frac{\chi(n)}{n^{\alpha}}$ est convergente alors, d'après le \textsc{Lemme} \ref{L1}, la fonction $\rho(x) = f(e^x)e^{-\alpha x}$ est un élément de $\vb^*$. Or d'après le \textsc{Corollaire} \ref{C1} on déduit que $$\rho(x) \underset{x \to + \infty}{\sim} Res(\L_{\rho},\beta) e^{\beta x} .$$

Ce qui est $$f(e^x) \underset{x \to + \infty}{\sim}Res(\L_{\rho},\beta) e^{(\alpha + \beta)x} .$$

D'où $$f(x) \underset{x \to + \infty}{\sim} Res(\L_{\rho},\beta) x^{\alpha + \beta} .$$

Ce qu'il fallait démontrer.  }

 On rappelle que la fonction $\Lambda$ de \textbf{Von Mangoldt} est une fonction arithmétique définie sur $\N^*$ par $$\Lambda(n) := \begin{cases} \log p \quad \text{si} \ n=p^k , \quad k \in \N^*,p\in \P \\ \\ \quad 0 \qquad  \text{sinon} \end{cases}.$$
 
 La fonction définie pour tout $x \in  [1,+ \infty)$ tel que $x \neq p^k$ où $k\in \N^*$ et $p\in \P$ par $ \displaystyle \psi(x) = \sum_{n < x} \Lambda(n)$ est dite la fonction de \textbf{Chebyshev}, ainsi pour démontrer le théorème des nombres premiers il faut et il suffit de démontrer que $$\psi(x) \underset{x\to + \infty}{\sim} x .$$
 
 Il existe une forte relation entre la fonction $\zeta$ de \textbf{Riemann} et la fonction $\psi$, en effet $$-\frac{\zeta'(s)}{\zeta(s)} = \sum_{n \geq 1} \frac{\Lambda(n)}{n^s} = s \int_0^{+ \infty} \psi(e^x)e^{-sx}dx, \qquad \forall \sigma > 1.$$
 
 On rappelle aussi que la fonction $\zeta$ est holomorphe sur $\{\sigma \geq 1\}$ sauf au $s=1$ qui est le seul pôle simple de la fonction $\zeta$, ainsi d'après \textbf{Hadamard} et \textbf{De La Vallée Poussin} la fonction $\zeta$ ne s'annulle en aucun point du demi-plan $\{\sigma \geq 1 \}$. Alors on déduit que la fonction $-\frac{\zeta'}{\zeta}$ est holomorphe sur $\{\sigma \geq 1 \}$ sauf au point $s=1$ qui est le seul pôle simple de résidu égal à $1$. 
 
 Et on a le \textbf{Théorème des Nombres Premiers}:

\corol{$$\psi(x) \underset{x\to + \infty}{\sim} x .$$}{Soit $\alpha > 1$ un réel donné, on pose $$\rho(x) = \psi(e^x)e^{-\alpha x}, \qquad \forall x \in \R^+.$$

Alors puisque la fonction $- \frac{\zeta'}{\zeta}$ est holomorphe sur $\{\sigma > 1\}$ alors la série su terme général $\frac{\Lambda(n)}{n^{\alpha}}$ est convergente pour tout $\alpha > 1$.

D'une autre part, soit $\sigma > 1$ alors 

\begin{align*}\L_{\rho}(s) &= \int_0^{+ \infty}\rho(x)e^{-sx}dx \\ &= \int_0^{+\infty} \psi(e^x)e^{-\alpha x} e^{-sx}dx \\ &= \int_0^{+ \infty} \psi(e^x)e^{-(s+\alpha)x}dx \\&= - \frac{\zeta'(s+\alpha)}{(s+\alpha)\zeta(s+\alpha)} \end{align*}

Alors puisque la fonction $s\mapsto - \frac{\zeta'(s+\alpha)}{(s+\alpha)\zeta(s+\alpha)}$ est holomorphe sur $\{\sigma \geq 1-\alpha\}$ sauf au point $s=1-\alpha$ qui est le seul pôle simple de cette fonction où $$Res\left(- \frac{\zeta'(s+\alpha)}{(s+\alpha)\zeta(s+\alpha)},1-\alpha \right)=1 $$

Alors d'après le \textsc{Théorème} \ref{T2} on a $$\psi(x) \underset{x \to + \infty}{\sim}x^{\alpha + 1 - \alpha} $$

C'est à dire $$\psi(x) \underset{x \to + \infty}{\sim} x . $$

Ce qu'il fallait démontrer.
}


\begin{thebibliography}{plain}
  
   \bibitem{[CarBrun]}
          M.Carter et B. Van Brunt, \textit{The Lebesgue-Stieltjes Integral} a practical introduction, Springer (2000). 
   \bibitem{[Titch]}
          E.C. Titchmarsh,
          \emph{The Theory of The Riemann Zeta-Function} 2nd ed, revised by D. R. Heath-Brown, Oxford University Press (1986).
   \bibitem{[Rud]}
          Walter Rudin,
          \emph{Analyse réelle et complexe}, Troisième tirage MASSON Paris New York Barcelone Milan 1980. 

\end{thebibliography}
\end{document}